\documentclass{amsart}
\usepackage{amssymb,graphicx,cite,upref}
\usepackage[line,ps,color,poly,all]{xy}
\newxyColor{lightgrey}{.75}{gray}{}
\CompileMatrices
\newcommand{\Z}{\mathbb Z}%
\newcommand{\LL}{\mathcal{L}}%
\newcommand{\K}{\mathcal{K}}%
\DeclareMathOperator{\SU}{SU}\DeclareMathOperator{\id}{id}%
\DeclareMathOperator{\SL}{SL}%
\DeclareMathOperator{\Sp}{Sp}%
\DeclareMathOperator{\GL}{GL}%
\DeclareMathOperator{\OO}{O}\DeclareMathOperator{\vcd}{vcd}%
\newcommand{\isomto}{\overset{\sim}{\rightarrow}}%
\newcommand{\R}{\mathbb R}%
\newcommand{\HH}{\mathfrak{H}}%
\newcommand{\Q}{\mathbb Q}%
\newcommand{\C}{\mathbb C}%
\newcommand{\E}{\mathbb E}%
\newcommand{\G}{{G}}%
\newcommand{\lag}[1]{#1}%
\newcommand{\lsp}[1]{{}^{#1}\!}%
\DeclareMathOperator{\rank}{rank}%
\DeclareMathOperator{\stab}{Stab}%
\newcommand{\Stab}[1]{\stab_{\Gamma}\left(#1\right)}%
\newcommand{\I} {{\mathcal I}}%
\newcommand{\Par} {{\mathcal P}}%
\newcommand{\J} {{\mathcal J}}%
\renewcommand{\Im}{\operatorname{Im}}%
\newcommand{\ip}[2]{\left \langle #1,#2 \right \rangle}%
\newcommand{\vect}[3]{\begin{pmatrix} #1 \\ #2 \\ #3 \end{pmatrix}}%
\newcommand{\ve}{\vect{1}{0}{0}}%
\newcommand{\vw}{\vect{0}{0}{1}}%
\newcommand{\vsigma}{\vect{i}{1+i}{1}}%
\newcommand{\comment}[1]{}%
\newcommand{\temp}[1]{}%
\newcommand{\Label}[1]{{\label{#1} \comment{#1}}}%
\newcommand{\group}[1]{\mathbf{#1}}%
\DeclareMathOperator{\pr}{pr}%
\theoremstyle{plain}
\newtheorem{thm}{Theorem}[section]

\newtheorem*{mr}{Main Result}
\newtheorem{lem}[thm]{Lemma}
\newtheorem{prop}[thm]{Proposition}

\theoremstyle{definition}
\newtheorem{defn}[thm]{Definition}

\newcommand{\OOO}{\mathcal{O}}
\begin{document}


\title{On the existence of spines for $\Q$-rank 1 groups}
\author{Dan Yasaki}
\address{Department of Mathematics and Statistics\\Lederle Graduate Research Tower\\ University of Massachusetts\\Amherst, MA 01003-9305}
\email{yasaki@math.umass.edu}
\date{}
\thanks{The original manuscript was prepared with the \AmS-\LaTeX\ macro
system and the \Xy-pic\ package.}
\begin{abstract}
Let $X=\Gamma \backslash G /K$ be an arithmetic quotient of a symmetric space of non-compact type.  In the case that $G$ has $\Q$-rank 1, we construct $\Gamma$-equivariant deformation retractions of $D=G/K$ onto a set $D_0$.  We prove that $D_0$ is a spine, having dimension equal to the virtual cohomological dimension of $\Gamma$.  In fact, there is a $(k-1)$-parameter family of such deformations retractions, where $k$ is the number of $\Gamma$-conjugacy classes of rational parabolic subgroups of $G$.  The construction of the spine also gives a way to construct an exact fundamental domain for $\Gamma$.  
\end{abstract}
\maketitle

\bibliographystyle{../amsplain_initials}

\begin{section}{Introduction}\Label{chapter:introduction}
Let $D=G/K$ be a symmetric space of non-compact type, where $G$ is the group of  real points of an semisimple algebraic group $\group{G}$ defined over $\Q$.  Let $\Gamma$ be an arithmetic subgroup of the rational points $\group{G}(\Q)$.  Let $(E,\rho)$ be a $\Gamma$-module over $R$.  If $\Gamma$ is torsion-free, the locally symmetric space $\Gamma \backslash D$ is a $K(\Gamma,1)$ since $D$ is contractible, and the group cohomology of $\Gamma$ is isomorphic to the cohomology of the locally symmetric space, i.e. $H^*(\Gamma,E)\cong H^*(\Gamma \backslash D; \E)$, where $\E$ denotes the local system defined by $(E,\rho)$ on $\Gamma \backslash D$.  When $\Gamma$ has torsion, the correct treatment involves the language of orbifolds, but the isomorphism of cohomology is still valid by using a suitable sheaf $\E$ as long as the orders of the torsion elements of $\Gamma$ are invertible in $R$. 

The \emph{virtual cohomological dimension} ($\vcd$) of $\group{G}$ is the smallest integer $p$ such that cohomology of $\Gamma \backslash D$ vanishes in degrees above $p$, where $\Gamma \subset \group{G}(\Q) $ is any torsion-free arithmetic subgroup.  Borel and Serre \cite{BS} show that the discrepancy between the dimension of $D$ and the $\vcd(\group{G})$ is given by the $\Q$-rank of $\group{G}$, the dimension of a maximal $\Q$-split torus in $\group{G}$.  Thus one can hope to find a $\Gamma$-equivariant deformation retract $D_0 \subset D$ of dimension equal to the virtual cohomological dimension of $\group{G}$.  When such a subset exists, it is called a \emph{spine}.

Spines have been constructed for many groups \cite{Sou,Men,Vog,Brown,A,MM, Bat,LSz}.  In \cite{A2}, Ash describes the \emph{well-rounded retract}, a method for constructing a spine for all linear symmetric spaces.  Ash and McConnell extend \cite{A2} to the Borel-Serre compactification in \cite{AM}.  The well-rounded retract works for algebraic groups $\group{G}$ where the real points are isomorphic to a product of the following groups \cite{FK}:
\begin{enumerate}
\item $\GL_n(\R)$.
\item $\GL_n(\C)$.
\item $\GL_n(\mathbb{H})$.
\item $\OO(1,n-1) \times \R^\times$.
\item The non-compact Lie Group with Lie algebra $\mathfrak{e}_{6(-26)}\oplus \R$.
\end{enumerate}
The retract has been used in the computation of cohomology \cite{Sou,LSz,A80,Men,Vog,SVog,AGG,AM1,AM2,vGT}.

The well-rounded retract proves the existence and gives a method of explicitly describing spines in linear symmetric spaces.  However, for non-linear symmetric spaces, no general technique to construct spines is known.  In fact, there were no examples until MacPherson and McConnell  \cite{MM} constructed a spine in the Siegel upper half-space for the $\Q$-rank 2 group $\Sp_4(\R)$.  
\begin{figure}
\begin{center}
\includegraphics[height=0.333\textwidth,width=\textwidth]{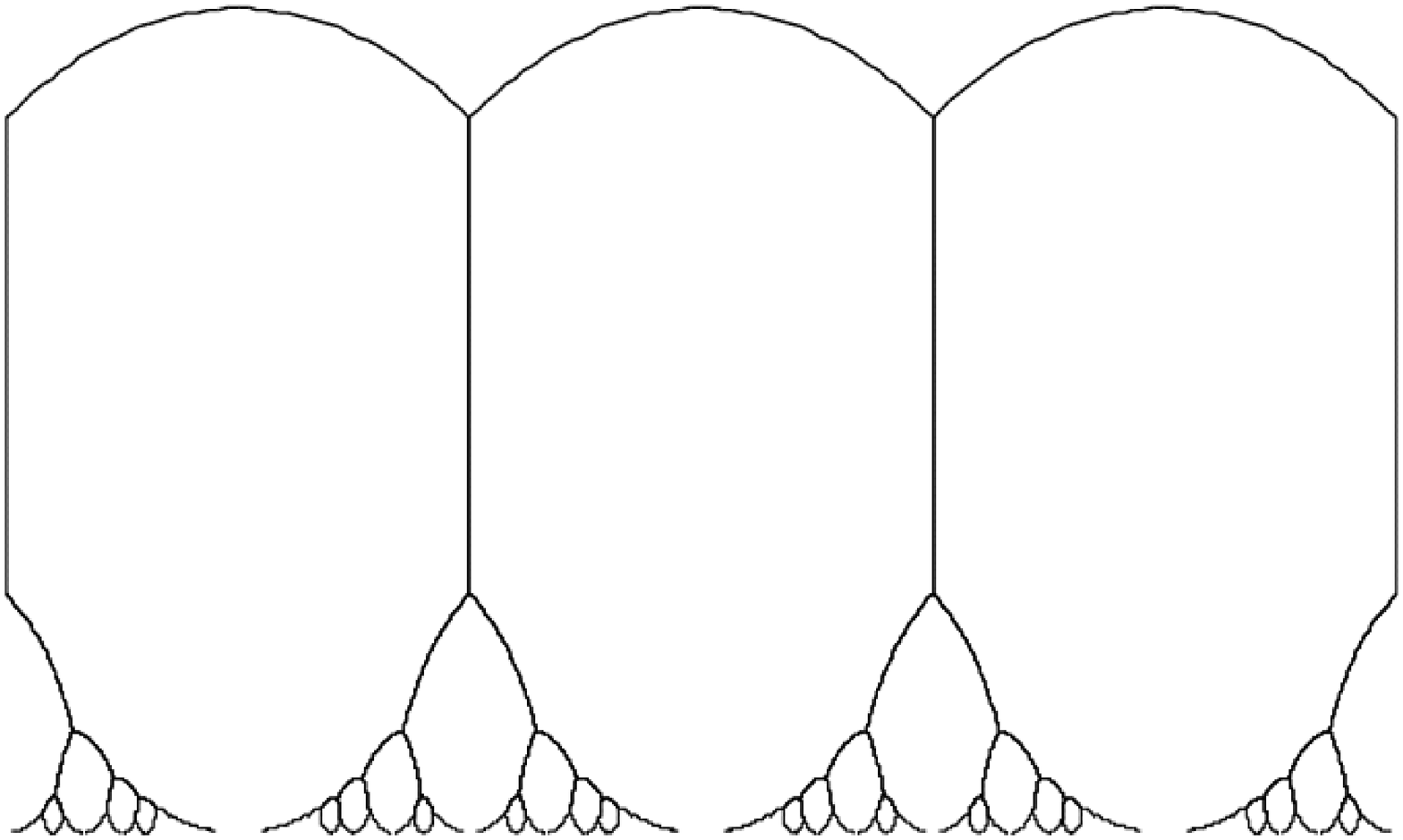}
\caption{Spine for $\SL_2(\R)$}
\Label{fig:sl2spine}
\end{center}
\end{figure}

In this paper, we deal with the case when $\group{G}$ has $\Q$-rank 1, and use a family of exhaustion functions to define a spine.  We use the exhaustion functions to construct a deformation retraction of $D$ onto a spine.  Each exhaustion function can be thought of as a measure of height with respect to a rational parabolic subgroup of $\group{G}$.  The existence of such functions is not new.  Siegel defined a notion of distance from a cusp for $\SL_2(\OOO_k)$ \cite{Sie}.  One can show that his distance functions are a power of our exhaustion functions.  More directly, the exhaustion functions come out of Saper's work on tilings in \cite{Sa}.  In fact, our exhaustion functions are nothing more than the composition of his \emph{normalized parameters} (in the $\Q$-rank 1 case) with the rational root.  We use the exhaustion functions to construct a deformation retraction of $D$ onto a spine.  More precisely,
\begin{mr}
Let $D=G/K$ be a symmetric space of non-compact type, where $G$ is the group of  real points of a semisimple algebraic group $\group{G}$ defined over $\Q$ with $\Q$-$\rank$ $1$.  Let $\Gamma$ be an arithmetic subgroup of the rational points $\group{G}(\Q)$.  There exists a $\Gamma$-equivariant retract of the symmetric space onto a set $D_0\subset D$ with the following properties:
\begin{enumerate}
\item $D_0$ is a locally finite union of semi-analytic sets.
\item $\dim(\Gamma \backslash D_0)=\vcd(\Gamma)$.
\item $D_0$ has a decomposition $D_0=\coprod D'(\I)$, where $\I$ ranges over certain subsets \textup{(}of order at least $2$\textup{)} of parabolic $\Q$-subgroups.
\item The decomposition satisfies $\gamma \cdot D'(\I)=D'(\lsp{\gamma} \I)$ for every $\gamma \in \Gamma$.  
\item The quotient $\Gamma \backslash D_0$ is compact.
\end{enumerate}
In fact, there is a $(k-1)$-parameter family of different retractions, where $k$ is the number of $\Gamma$-conjugacy classes of parabolic $\Q$-subgroups.  
\end{mr}

For $\group{G}=\SL_2$, our technique yields the same deformation retraction as the well-rounded retract, which is an infinite trivalent tree in the Poincar\'e upper half-plane (Figure~\ref{fig:sl2spine}).  In this case, the exhaustion function corresponding to the cusp $i\infty$ is simply $f_{i\infty}(z)=\Im(z)$.  Further comparisons with other known results are given in Section~\ref{sec:examples}.

Since each $D'(\I)$ is a semi-analytic set, it follows from \cite{Loj} that the decomposition of the spine $D_0=\coprod D'(\I)$ may be refined to a regular cell complex.  Thus, the cohomology can be described by finite combinatorial data.  In \cite{Yaspicard} we develop machinery for the computation of cohomology using $D_0$, and use it together with the results of this paper to investigate the cohomology of $\SU(2,1)$ over the Gaussian integers.

Sections~\ref{sec:notationandbackground} and \ref{sec:definitions} set notation and define the main objects of this paper.  Sections~\ref{sec:rep} and \ref{sec:properties} give another interpretation of the exhaustion functions and prove some of their properties.  The main results are presented in Section~\ref{sec:spine}.  Section~\ref{sec:separated} introduces the notion of a strictly separated linear algebraic group, and show that for these groups, each $D'(\I)$ is a smooth, contractible submanifold.  Section~\ref{sec:examples} concludes by looking at a few examples of spines in low dimensional cases.

I would like to thank my thesis advisor, Les Saper, for his insight into this work.  I would also like to thank Paul Gunnells for helpful conversations. 
\end{section}
\begin{section}{Notation and background}\label{sec:notationandbackground}
In order to set notation, we briefly recall without proof some standard results regarding algebraic groups over $\Q$, the geodesic action, and the Borel-Serre compactification \cite{BS}.  We follow the exposition of \cite{Sa}.  Throughout this paper, $G$ is the identity component of the real points of a semisimple $\Q$-rank $1$ algebraic group defined over $\Q$.  In order to lighten the notation and exposition, we will notate the algebraic group and its group of real points by the same Roman type $(G=G(\R))$, and we may refer to the algebraic group when properly we should referring to the group of real points.  For example, we will speak of parabolic $\Q$-subgroups of $G$, when we properly should be referring to the group of real points of a parabolic $\Q$-subgroup of the algebraic group $G$.
\begin{subsection}{Algebraic groups defined over $\Q$}\Label{sec:background}
For a reductive algebraic group $H$ which is defined over $\Q$, let ${S_H}$ denote the maximal $\Q$-split torus in the center of $H$, and set $A_H={S_H}(\R)^0$.  Set 
\begin{equation*}
\lsp{0}H=\bigcap_{\chi \in X(H)_\Q} \ker(\chi^2),
\end{equation*}
where $X(H)_\Q$ denotes the rational characters of $H$ defined over $\Q$.  Then $H$ splits as a direct product
\[
H=A_H \times \lsp{0}H.
\]
The group $\lsp{0}H \subset H$ contains all compact and arithmetic subgroups of $H$.

For a parabolic subgroup $\lag{P}\subset \G$, let $\lag{N_P}$ denote its unipotent radical of and let 
\begin{equation*}
\nu_P:\lag{P} \to \lag{L_P}=\lag{P}/\lag{N_P}
\end{equation*}
denote the projection to the Levi quotient.  Let $M_P$ denote the group $\lsp{0}L_P$.  If $\lag{S_P}\subset \lag{L_P}$ denotes the maximal $\Q$-split torus in the center of $\lag{L_P}$, then $L_P$ splits as a commuting direct product $L_P=A_PM_P$, where $A_P=\lag{S_P}(\R)^0$.  Note that $L_P$ is the centralizer of $A_P$, and the connected component of the center of $M_P$ is a $\Q$-anisotropic torus.   

Let ${}_\Q \Delta _P =\{\alpha_P\}$ denote the simple root of the adjoint action of $A_P$ on $\mathfrak{n}_P$, the Lie algebra of $N_P$.  The roots will be viewed as characters of $A_P$ and as elements of $\mathfrak{a}_P^*$.  Since $G$ has $\Q$-rank 1, $\alpha_P$ gives an isomorphism $A_P \to \R_{>0}$.  Let 
\begin{equation}\Label{eq:Raction}
\psi_P:\R_{>0}\to A_P
\end{equation}
be the isomorphism of groups given by $\psi_P(t)=\alpha_P^{-1}(t)$.

For two rational parabolic subgroups $P$ and $Q$, there is a canonical isomorphism 
\begin{equation*}
A_P \isomto A_Q
\end{equation*}
induced by an element of $G(\Q)$.  In particular, the characters $\{\alpha_P\}_{P\in \Par}$ can be identified and will be denoted $\alpha$.

Any lift $i:\lag{L_P}\to \lag{P}$ determines a \emph{Langlands decomposition}, which is a semi-direct product
\begin{equation*}
P=N_Pi(A_PM_P).
\end{equation*}
Note that $L_P$ is the centralizer of $A_P$ and the connected component of the center of $M_P$ is a $\Q$-anisotropic torus.   

Let $K \subset G$ be a maximal compact subgroup and define $D=G/K$.  There is a unique basepoint $x\in D$ such that $K=\stab_G(x)$.  This choice of $K$ also determines the following data \cite{GoTai}:
\begin{enumerate}
\item A maximal compact subgroup $K_P = K \cap P \subset G$ and a diffeomorphism $P/K_P \to D$.
\item A Cartan involution $\theta_x: G \to G$ such that $K$ is the subgroup fixed by $\theta_x$.
\item A unique $\theta_x$-equivariant lifting $i_{x}: \lag{L_P} \to \lag{P}$.  For a subset $T \subset L_P$, denote its lift by $T(x)$.
\end{enumerate}

\begin{defn}\cite{GoHMac}\label{defn:rationalbasepoint}
Let $P \subset G$ be a rational parabolic subgroup and $x$ a point in $D$.  If the lift $L_P(x)$ is an algebraic subgroup of $P$, then $x$ is a \emph{rational basepoint for $P$}.  
\end{defn}
If the basepoint $x$ can be chosen so that the associated maximal compact subgroup $K$ is defined over $\Q$, then $x$ is rational for all rational parabolic subgroups of $G$.
\begin{defn}\Label{defn:QWeyl}
Let $S$ be a maximal $\Q$-split torus of $G$.  The \emph{Weyl group over $\Q$} or \emph{$\Q$-Weyl group} is the quotient of the normalizer $\mathcal{N}(S)$ of $S$ by the centralizer $\mathcal{Z}(S)$ of $S$, and is denoted
\begin{equation*}
{}_\Q W=\mathcal{N}(S)/\mathcal{Z}(S).
\end{equation*}
\end{defn}

The following is the rational version of the standard \emph{Bruhat Decomposition} for Lie groups. 
\begin{thm}\textup{\cite{Bo}} \label{thm:Bruhat}
Let $P\subset G$ be a minimal rational parabolic subgroup.  Then $G(\Q)$ is the disjoint union of the classes $P(\Q)wP(\Q)$ with $[w]\in {}_\Q W$.  In particular, given $g \in G(\Q)$, there exists $u_g \in N_P$, $[w] \in {}_\Q W$, and $p_g \in P$ such that $g= u_g w p_g$.
\end{thm}
\end{subsection}

\begin{subsection}{Borel-Serre compactification}\cite{BS}
Let $x$ denote the point of $D$ fixed by $K$.  Then $P \in \Par$ acts transitively on D, so every point $z \in D$ can be written as $z=p \cdot x$, for some $p \in P$.  The geodesic action of $A_P$ on $D$ is given by
\begin{equation*}
a \circ z = (p \tilde a) \cdot x,
\end{equation*}
where $\tilde a$ is the image of $a$ in $A_P(x)$.  The geodesic action commutes with the usual action of $P$ on $D$ and is independent of the choice of basepoint.  Let $A_P \times \lsp{0}P$ act on $D$ by $(a,p) \cdot z = p \cdot (a \circ z)=a \circ (p \cdot z)$.  Then there is an analytic isomorphism of $A_P \times \lsp{0}P$-homogeneous spaces 
\begin{equation}\Label{eqn:D(P)}
(a_{P,x},q_P):D \isomto A_P \times e(P) ,\end{equation}
where $e(P)=A_P \backslash D$ is the quotient of $D$ by the geodesic action of $A_P.$  Normalize $a_{P,x}$ so that $a_{P,x}(x)=e$.  Via~\eqref{eqn:D(P)}, $D$ is a trivial principal $A_P$-bundle with \emph{canonical cross-sections} given by the orbits of $\lsp{0} P$.  

The Borel-Serre compactification is then constructed as follows:  The simple root $\alpha_P$ induces an isomorphism $A_P \isomto \R_{>0}$ defined by $a \mapsto a^\alpha$.  The \emph{partial bordification} $D^c(P)$ associated to $P$ is defined to be $\overline{A}_P \times_{A_P} D$.  Equivalently, extend \eqref{eqn:D(P)} to 
\begin{equation*}
(a_{P,x},q_P):D^c(P) \isomto \overline{A}_P \times e(P).
\end{equation*} 
The Borel-Serre compactification $\overline{D} \equiv \bigcup_{P \in \Par} D^c(P)$ is then given the unique structure of an analytic manifold with boundary so that each $D^c(P)$ is an open submanifold with boundary.

The action of $\G(\Q)$ on $D$ extends to an action on $\overline{D}$, and for an arithmetic subgroup $\Gamma \subset \G(\Q)$, the quotient $\Gamma \backslash \overline{D}$ is compact.
\end{subsection}

\begin{subsection}{$\Gamma$-conjugacy classes of parabolic $\Q$-subgroups}
Fix a proper parabolic $\Q$-subgroup $P$ of $\G$.  The parabolic $\Q$-subgroups of $G$ are all conjugate to ${P}$ via elements of $\G(\Q)$.   Thus, the set of $\Q$-parabolic subgroups $\Par$ is in one-to-one correspondence with points of $\G(\Q) / {P}(\Q)$, where $[g] \in \G(\Q) / {P}(\Q)$ corresponds to the parabolic $\Q$-subgroup $\lsp{g} P = P^{g^{-1}} \equiv gPg^{-1}$.  For an arithmetic subgroup $\Gamma \subset \G(\Q)$, the double coset space $\Gamma \backslash \G(\Q) / {P}(\Q)$ is finite \cite{BoHC} and the number of elements in the space is known as the \emph{class number}.  In particular, there are only finitely many $\Gamma$-conjugacy classes of parabolic $\Q$-subgroups of $G$.  For $P \in \Par$, let $\Gamma_{P} \equiv  \Gamma \cap P \subset \lsp{0} P$.  Then there is a bijection between $\Gamma/\Gamma_{P}$ and the parabolic subgroups that are $\Gamma$-conjugate to $P$ via $\gamma \Gamma_P \mapsto \lsp{\gamma} P$.
\end{subsection}
\end{section}

\begin{section}{Main definitions}\Label{sec:definitions}
In this section, we define a family $\{f_P\}_{P\in \Par}$ of \emph{exhaustion functions} depending on a \emph{$\Gamma$-invariant parameter}.  We also define subsets of $D$ associated to the family of exhaustion functions.  In Section~\ref{sec:spine}, these functions are used to give a $\Gamma$-equivariant retraction of $D$ onto a codimension 1 set $D_0$.
\begin{defn}
A \emph{parameter} is a family $\{O_P\}_{P\in \Par}$ of closed submanifolds of $D$ such that each $O_P$ has the form 
\begin{equation*}
O_P=\lsp{0} P \cdot x_P \quad \text{for some $x_P \in D$.}
\end{equation*}
A parameter is \emph{$\Gamma$-invariant} if 
\begin{equation*}
\gamma \cdot O_P=O_{\lsp{\gamma}P} \quad \text{for all $\gamma \in \Gamma$ and $P\in \Par$.}
\end{equation*}
\end{defn}
Recall that by conjugation by $\G(\Q)$, one can canonically identify the $A_P$ and corresponding simple root $\alpha_P$ for different $P\in \Par$.  Denote the simple root by $\alpha$ and view it as a character on each $A_P$.  

For $P \subset G$ a rational parabolic subgroup, let 
\begin{equation}
a_P:D \times D \to A_P
\end{equation}
be the map $a_P(z,x)=a_{P,x}(z)$. Thus $a_P(z,x)$ can be viewed as the amount of geodesic action required to push $\lsp{0}P \cdot x$ to $\lsp{0}P \cdot z$.  The following is immediate from the definitions and the $\lsp{0}P$-invariance of $a_{P,x}$.
\begin{prop} \textup{\cite{BS}}\Label{prop:aproperty}
Let $P\subset G$ be a rational parabolic subgroup.  Then 
\begin{enumerate}
\item $a_P(z,x)=a_P(x,z)^{-1} \quad \text{for all $z, x \in D$.}$
\item $a_{\lsp{g}P}(g \cdot z, g\cdot x)=a_P(z,x) \quad \text{for all $z, \ x \in D$ and $g \in \G(\Q)$.}$
\item $a_P(z,x)=a_P(p \cdot z,x)=a_P(z,p \cdot x)$ \quad for all $p \in \lsp{0}P$.
\item $a_P(z,s)a_P(s,x)=a_P(z,x)$ \quad for all $z,\ x,$ and $s$ in $D$.
\end{enumerate}
\end{prop}

Given a point $x\in D$ and rational parabolic subgroup $P \in \Par$, consider the function $f:D \to \R_{>0}$ given by $f(z)=a_P(z,x)^\alpha$.  Proposition~\ref{prop:aproperty} implies that $f$ only depends on the orbit $\lsp{0}P\cdot x$.  This leads to the following definitions.
\begin{defn}\Label{defn:exhaustion}
The \emph{exhaustion functions associated a parameter $\{\lsp{0}P\cdot x_P\}_{P\in \Par}$} is the family $\{f_P\}_{P\in \Par}$ of functions $f_P:D \to \R_{>0}$ given by
\begin{equation*}
f_P(z)=a_P(z,x_P)^\alpha.
\end{equation*}
A family of exhaustion functions associated to a $\Gamma$-invariant parameter is called a family of \emph{$\Gamma$-invariant exhaustion functions}.
\end{defn}
Notice that Proposition~\ref{prop:aproperty} implies that 
\begin{equation}\Label{eq:finvariant}
f_{\lsp{g}P}(g \cdot z)=f_P(z)f_{\lsp{g}P}(g \cdot x_P) \quad \text{for all $g \in G(\Q)$,}
\end{equation}
so that in particular, for exhaustion functions associated to a $\Gamma$-invariant parameter,
\begin{equation}\Label{eq:gammainvariant}
f_{\lsp{\gamma}P}(\gamma \cdot z)=f_P(z) \quad \text{for all $\gamma \in \Gamma$.}
\end{equation}

For a parabolic $P$, define $D(P)\subset D$ to be the set of $z \in D$ such that $ f_P(z) \geq  f_Q(z)$ for every $Q \in \Par \setminus \{P\}$.  More generally, for a subset  $\I \subseteq \Par$, 
\begin{align}
E(\I)&=\{z\in D\; |\;  f_P(z)= f_Q(z) \text{ for every pair } P,Q \ \in \I \}\\
D(\I)&=\bigcap_{P\in \I} D(P)\\ 
D'(\I)&=D(\I) \setminus \bigcup_{\I' \supsetneq \I} D(\I').
\end{align}
It follows that $D'(\I)\subseteq D(\I) \subset E(\I) \text{ and } D(\I)=\coprod_{\tilde{\I} \supseteq \I}D'(\tilde{\I})$.  

\begin{defn}\label{defn:degenerate}
The set $D'(\I)$ will be called a \emph{degenerate tile}.  Let $f_\I$ denote the restriction to $E(\I)$ of $f_P$ for $P\in \I$. 
\end{defn}

\begin{defn}
Let $\I \subseteq \Par$, $P \in \I$, and $z \in E(\I)$.  Then $z$ is called a {\it first contact for $\I$} if $f_\I(z)$ is a global maximum of $f_\I$ on $E(\I)$.
\end{defn}

\begin{defn}
A subset $\I \subset \Par$ is called {\it admissible} if $D(\I)$ is non-empty and {\it strongly admissible} if $D'(\I)$ is non-empty.
\end{defn}

\begin{prop}\Label{prop:tiling}
Let $\mathcal S$ denote the collection of strongly admissible subsets of $\Par$.  Then the symmetric space has a $\Gamma$-invariant degenerate tiling 
\[D = \coprod_{\I \in {\mathcal S}} D'(\I),\]
such that $\gamma \cdot D'(\I)=D'(\lsp{\gamma} \I)$ for all $\gamma \in \Gamma$ and $\I\in \mathcal{S}$.
\end{prop}
\begin{defn}\Label{defn:spine}
Given a family of $\Gamma$-invariant exhaustion functions, define a subset $D_0 \subset D$ by
\begin{equation*}
D_0=\coprod_{\substack{\I \in {\mathcal S}\\ |\I| >1}} D'(\I).
\end{equation*} 
\end{defn}
Let $f_{D_0}$ denote the function on $D_0$ given by
\begin{equation}\Label{eq:fD0}
f_{D_0}(z)=f_\I(z) \quad \text{for $z\in D(\I)$.}
\end{equation}
\end{section}
\begin{section}{Exhaustion functions via representation theory}\Label{sec:rep}
In this section we describe a systematic way to construct exhaustion functions. 

Fix a parabolic $\Q$-subgroup $P \subset G$.  Choose a rational basepoint ${x^*}\in D$ for $P$ so that the $\theta_{x^*}$-stable lift $A_P({x^*})$ of $A_P$ to $P$ is the connected component of the real points of a maximally $\Q$-split torus of $\G$.  Let $M_P({x^*})$ denote the $\theta_{x^*}$-stable lift of $M_P=\lsp{0}L_P$ to $P$, and let $\mathfrak{a}_P$ denote the Lie algebra of $A_P({x^*})$.  The simple $\Q$-root $\alpha$ can be viewed as a linear functional on $\mathfrak{a}_P$.  Let $\mathfrak{g}_\C$ denote the complexification of $\mathfrak{g}$.  Fix a Cartan subalgebra containing the complexification of $\mathfrak{a}_P$, and let $\{\tilde \alpha_1,\ldots \tilde \alpha_m\}$ be the simple roots.  Order the simple roots so that $\tilde \alpha_j|_{\mathfrak{a}_P}=\alpha$ for $1 \leq j \leq l$ and  $\tilde \alpha_j|_{\mathfrak{a}_P}=0$ for $l+1 \leq j \leq m$.     
Let $\{\tilde \omega_1,\ldots,\tilde \omega_m\}$ be the associated fundamental weights.  Consider the weight $\tilde \omega=s\ \sum_{j=1}^l \tilde \omega_j$.  Then for $s$ a sufficiently large positive integer, $\tilde\omega$ is the highest weight of a finite-dimensional (strongly rational) representation of $\G$ \cite{BT1}. 

Let $(V,\pi)$ denote this representation, and $\omega$ the restriction of $\tilde \omega$ to $\mathfrak{a}_P$.  The restriction to $\mathfrak{a}_P$ of weights are called \emph{restricted $\Q$-weights}.  Every restricted $\Q$-weight of $V$ is of the form $\omega-j \alpha$, where $j$ is a non-negative integer and the lowest restricted $\Q$-weight is $-\omega$.  Thus $V$ has a decomposition as a direct sum of restricted $\Q$-weight spaces,
\begin{equation} \Label{eq:Vdecomp}
V=\bigoplus_\lambda V_\lambda = \bigoplus_{k=0}^N V_{\omega -k \alpha}, \quad \text{where $\omega-N\alpha=-\omega$.}
\end{equation} 
Fix an \emph{admissible inner product} $\ip{\cdot}{\cdot}$ on $V$ \cite{BoWa}, that is, one for which 
\begin{equation*}
\pi(g)^*=\pi(\theta_{x^*} g)^{-1} \quad  \text{for all $g \in G$,}
\end{equation*}
so that in particular, 
\begin{equation*}
\pi(k)^*=\pi(k)^{-1} \quad  \text{for all $k \in K$.}
\end{equation*}

Since the $\Q$-Weyl group ${}_\Q W$ is finite, by averaging over ${}_\Q W$ one can arrange that the inner product is also ${}_\Q W$-invariant.  The data $(V,\pi,\ip{\cdot}{\cdot},x^*)$ will be called a \emph{$P$-adapted representation}.
\begin{prop}\Label{prop:Padaptedrep}
Let $P\subset \G$ be a parabolic $\Q$-subgroup with rational basepoint ${x^*}$, and let $(V,\pi,\ip{\cdot}{\cdot},x^*)$ be a $P$-adapted representation. Then $M_P({x^*})$ preserves the restricted $\Q$-weight spaces of $V$.  The highest and lowest restricted weight spaces are one-dimensional, and $M_P({x^*})$ preserves length on both.  
\end{prop}
\begin{proof}
Since $M_P({x^*})$ centralizes $A_P({x^*})$, it follows that $M_P({x^*})$ preserves the restricted weight space decomposition of $V$.  Since $\tilde \omega$ is orthogonal to $\tilde \alpha_j$ for $l+1 \leq j\leq m$, it follows that the highest weight space $V_\omega$ is $1$-dimensional.  Furthermore, this implies that the connected component of $M_P({x^*})$ acts trivially on $V_\omega$.  $M_P(x^*)$ only has finitely many connected components because $M_P \times A_P \cong L_P$, the real points of the Levi quotient, which has only finitely many connected components since it is Zariski connected \cite{Wh}.  It follows that $M_P(x^*)$ preserves length when restricted to $V_\omega$.  An analogous argument shows that $V_{-\omega}$ is $1$-dimensional, and $M_P(x^*)$ preserves length on $V_{-\omega}$.
\end{proof}

Write $z\in D$ as $z= p \cdot x$ with $p\in P$.  Using Langlands decomposition, write $p$ as $u \tilde a_P(z,{x^*}) m$, where $u \in N_P$, $\tilde a_P(z,{x^*})\in A_P({x^*})$, and $m \in M_P({x^*})$.  Fix $g \in \G(\Q)$.  From \eqref{eq:Vdecomp} there exists unit vectors $v_k \in V_{\omega -k \alpha}$ and constants $c_k(z) \in \C$ such that 
\begin{equation}
\pi(p^{-1} g) v=\pi(\tilde a_P(z,{x^*}))^{-1}\sum_k c_k(z) v_k.
\end{equation}  
Note that $c_k(z)$ is only defined up to multiplication by norm 1 scalars.  Since $M_P({x^*})$ centralizes $A_P({x^*})$, the $c_k(z)$ only depend on $m$ and $u$.  Thus for each fixed $g$, the norm $|c_k|$ can be viewed as a function from $D \to \R_{\geq 0}$ that is invariant under the geodesic action of $A_P$. 

\begin{prop} \Label{prop:fPexplicit}
Fix $g \in G(\Q)$ and two distinct rational parabolic subgroups $P$ and $Q=\lsp{g} P$.  Fix a rational basepoint ${x^*}$ for $P$ and let $(V,\pi,\ip{\cdot}{\cdot}, x^*)$ be a $P$-adapted representation.  Then the functions $|c_k|:D \to \R$ \textup(defined above\textup) satisfy 
\begin{equation*}
f_Q(z)=\frac{f_Q(g \cdot x_P)f_P(z)}{\left(\sum_{k=0}^N |c_k(z)|^2 f_P({x^*})^{-2k}f_P(z)^{2k}\right)^{1/N}} \quad \text{for all $z \in D$.}
\end{equation*}
\end{prop}
\begin{proof}
From the definition of the $P$-adapted representation, $\alpha=\frac{2 \omega}{N}$.  Write $z \in D$ as $z=p \cdot {x^*}$, where $p\in P$.  Let $v\in V$ denote a norm 1 highest weight vector.  Since $N_P$ fixes $v$, Proposition~\ref{prop:Padaptedrep} implies 
\begin{equation}\Label{eq:pi}
\|\pi(p^{-1})v\|^{-2/N}=\|\pi(m^{-1}\tilde a_P(z,{x^*})^{-1}u^{-1})v\|^{-2/N}=a_P(z,{x^*})^\alpha.
\end{equation}
Then from the definition of $f_P$ and \eqref{eq:pi},
\begin{equation}\Label{eq:fPrep}
f_P(z)=a_P(z,{x^*})^\alpha f_P({x^*})=\|\pi(p^{-1})v\|^{-2/N} f_P({x^*}).
\end{equation}
Then \eqref{eq:finvariant}, \eqref{eq:fPrep}, and the orthogonality of the restricted weight spaces imply,
\begin{align*}
f_Q(z)&=
f_P(g^{-1}\cdot z)f_Q(g \cdot x_P)\\
&=f_Q(g \cdot x_P)f_P({x^*})\biggl(\biggl\|\pi(\tilde a_P(z,{x^*}))^{-1}\sum_{k=0}^N c_k(z)v_k\biggr\|^2\biggr)^{-1/N}\\
&=f_Q(g \cdot x_P)f_P({x^*})\biggl(\biggl\|\sum_{k=0}^Nc_k(z) a_P(z,{x^*})^{-\omega + k\alpha}v_k\biggr\|^2\biggr)^{-1/N}\\
&=f_Q(g \cdot x_P)f_P({x^*})\biggl(\sum_{k=0}^N|c_k(z)|^2a_P(z,{x^*})^{-2\omega +2k\alpha}\biggr)^{-1/N}\\
&=f_Q(g \cdot x_P)f_P({x^*})a_P(z,{x^*})^{2\omega/N}\biggl(\sum_{k=0}^N|c_k(z)|^2a_P(z,{x^*})^{2k\alpha}\biggr)^{-1/N}\\
&=f_Q(g \cdot x_P)f_P(z)\biggl(\sum_{k=0}^N |c_k(z)|^2 f_P({x^*})^{-2k}f_P(z)^{2k}\biggr)^{-1/N}.\qedhere
\end{align*}
\end{proof}

\begin{prop}\Label{prop:ckexplicit}
Let $P$ and $Q=\lsp{g}P$ be two distinct $\Q$-parabolic subgroups with $g \in G(\Q)$.  Let $g=u_g w p_g$ be the $\Q$-Bruhat decomposition of $g$ as in Proposition~\ref{thm:Bruhat}.  Use the Langlands decomposition to express $z \in D$ as $z=u \tilde a_P(z,{x^*}) m \cdot {x^*}$ for a rational basepoint $x^*$ for $P$, $u\in N_P$, $m\in M_P({x^*})$, and $\tilde a_P(z,{x^*}) \in A_P({x^*})$. Then 
\begin{equation*}
|c_N(z)|=a_P(p_g \cdot {x^*},{x^*})^\omega.
\end{equation*}
If $|c_k(z)|=0$ for $k=0,\ldots,N-1$, then 
\begin{enumerate}
\item $u=u_g$.
\item $z$ is a first contact for $\{P,Q\}$.
\item $f_P(z)=f_Q(g \cdot x_P)^{1/2}a_P(p_g \cdot x^*,x^*)^{-\alpha/2}f_P(x^*)$.
\end{enumerate}
\end{prop}
\begin{proof}
From the definition of $|c_k(z)|$,
\begin{equation}
\begin{split}
|c_k(z)|&=\|\pr_{V_{\omega-k\alpha}}(\pi(m^{-1}u^{-1}g) v)\|\\
&=\|\pr_{V_{\omega-k\alpha}}(\pi(m^{-1}u^{-1} u_g w p_g) v)\|\\\Label{eq:ck}
&=a_P(p_g\cdot {x^*},{x^*})^\omega\|\pr_{V_{\omega-k\alpha}}(\pi(m^{-1}u^{-1} u_g w) v)\| .
\end{split}
\end{equation}
Notice that $\pi(w) v \in V_{-\omega}$ and for $n\in N_P$, 
\begin{equation*}
\pi(nw)v=\pi(w)v+ \text{higher weight vectors}. 
\end{equation*}
Since $M_P({x^*})$ preserves weight spaces and preserves length on $V_{-\omega}$ by Proposition~\ref{prop:Padaptedrep},
\begin{equation}
|c_N(z)|=a_P(p_g\cdot {x^*},{x^*})^\omega\|\pr_{V_{-\omega}}(\pi(m^{-1}u^{-1} u_g w) v)\|=a_P(p_g\cdot {x^*},{x^*})^\omega .
\Label{eqn:cN}
\end{equation}
Now suppose $|c_k(z)|=0$ for $k=0,\ldots,N-1$.  Then from \eqref{eq:ck}, $\pi(m^{-1}u^{-1} u_g w) v \in V_{-\omega}$.  Since $\pi(m^{-1})$ preserves weight spaces by Proposition~\ref{prop:Padaptedrep}, $\pi(u^{-1}u_g)$ preserves $V_{-\omega}$.  It follows that $u=u_g$.

Note that since $z\in E(\{P,Q\})$,
\[
f_Q(g \cdot x_P)^N=\sum_{k=0}^N \frac{|c_k(z)|^2}{f_P(x^*)^{2k}}f_P(z)^{2k} \leq \frac{|c_N(z)|^2}{f_P(x^*)^{2N}}f_P(z)^{2N}.
\]
Since $|c_N(z)|\neq 0$ and is independent of $z$ from \eqref{eqn:cN}, 
\begin{equation}
f_P(z)\leq \frac{f_Q(g \cdot x_P)^{1/2}f_P(x^*)}{|c_N(z)|^{1/N}}\quad \text{on $E(\{P,Q\})$.}
\end{equation}
The bound is attained when $u(z)=u_g$ and
\[
f_P(z)=\frac{f_Q(g \cdot x_P)^{1/2}f_P(x^*)}{|c_N(z)|^{1/N}}=\frac{f_Q(g \cdot x_P)^{1/2}f_P(x^*)}{a_P(p_g \cdot x^*,x^*)^{\alpha/2}}\qedhere
\]
\end{proof}
\end{section}
\begin{section}{Properties of Exhaustion Functions}\Label{sec:properties}
Notice that for a family of $\Gamma$-invariant exhaustion functions, 
\begin{equation}
f_{\lsp{\gamma}P}(\gamma \cdot x_P)=1 \quad \text{for every $\gamma \in \Gamma$.}
\end{equation}
Proposition~\ref{prop:aproperty} implies the following.
\begin{prop} \Label{prop:fproperty}
Let $P\subset G$ be a rational parabolic subgroup and $\{f_P\}_{P\in \Par}$ a family of exhaustion functions associated to a $\Gamma$-invariant parameter $\{\lsp{0}P\cdot x_P\}_{P\in \Par}.$  Then 
\begin{enumerate}
\item $f_P(x_P)=1$.
\item $f_P(a \circ z)=a^\alpha f_P(z)$  \quad for all $z\in D$ and $a\in A_P$.
\item $f_P(z)=a_P(z,x)^\alpha f_P(x) \quad \text{for all $z,\ x \in D$.}$
\item $f_P( p \cdot z)=f_P(z)$ \quad for all $z\in D$ and $p\in \lsp{0}P$.
\item $f_{\lsp{g}P}(z)=f_P(g^{-1}z)f_{\lsp{g}P}(g \cdot x_P)$ \quad for all $z\in D$ and $g \in G(\Q)$.
\end{enumerate}
\end{prop}

The following proposition follows from reduction theory \cite{BoHC}.
\begin{prop}\Label{prop:minf}
Let $\{f_P\}_{P\in \Par}$ be a $\Gamma$-invariant family of exhaustion functions.  Then there exists a constant $C > 0 $, depending on the $\Gamma$-invariant parameter, such that
\begin{equation*}
\sup_{P\in \Par}f_P(z) \geq C \quad \text{for all $z\in D$.}
\end{equation*} 
\end{prop}
\begin{lem} \Label{lem:chainrule}
Let $M$ be a Riemannian manifold and $\psi$ an isometry of $M$.  Let $f$ and $h$ be smooth functions on $M$ such that $h(x)=f(\psi(x))$.  Then \[\nabla h(x)=(d\psi^{-1})_{\psi(x)}\nabla f(\psi(x))\quad \text{for all $x\in M$.}\]
\end{lem}

\begin{prop}\Label{prop:constantlevelset}
Let $x,\ x' \in D$.  If $f_P(x)=f_P(x')$ then $\|\nabla f_P(x)\|=\|\nabla f_P(x')\|$. 
\end{prop}
\begin{proof}
If $ f_P(x)= f_P(x')$, then there exists $p\in \lsp{0}P$ such that $p\cdot x=x'$.  Since $ f_P$ is $\lsp{0}P$-invariant, $ f_P(p\cdot z)= f_P(z)$ for all $z\in D$.  The result then follows from Lemma~\ref{lem:chainrule} by noting that $\lsp{0}P \subset G$ acts by isometries on $D$.
\end{proof}

\begin{lem}\Label{lem:geodesicisometry}
Let $P \subset G$ be a rational parabolic subgroup.   The geodesic action of $A_P$ acts as an isometry on the $1$-dimensional tangent space $T_z(A_P \circ z)$. 
\end{lem}
\begin{proof}
Since the geodesic action of $A_P$ is equal to the regular action of its $\theta_z$-invariant lift on the orbit $A_P \circ z$, the geodesic action acts as an isometry on the $1$-dimensional tangent space $T_{z}(A_P \circ z)$.  Furthermore, the geodesic action of $A_P$ commutes with the regular action of $P$ on $D$.  The result follows.
\end{proof}

\begin{prop}\Label{prop:norminvariant}
Fix two distinct rational parabolic subgroups $P$ and $Q$.  Then 
\begin{equation*}
f_P(z)\|\nabla f_Q(z)\|=f_Q(z)\|\nabla f_P(z)\| \quad \text{for all $z\in D$.}
\end{equation*} 
In particular, if $z\in E(\{P,Q\})$, then $\|\nabla  f_P(z)\|=\|\nabla  f_Q(z)\|$. 
\end{prop}
\begin{proof}
There exists a $g\in \G(\Q)$ such that $Q=\lsp{g}P$.  By Proposition~\ref{prop:fproperty},
\begin{equation} 
f_Q(z)=f_P(g^{-1}\cdot z)f_Q(g \cdot x_P).
\end{equation}
Then Lemma~\ref{lem:chainrule} implies that 
\begin{equation}\Label{eq:e1}
\nabla f_Q(z)=f_Q(g \cdot x_P)(dL_g)_{g^{-1}\cdot z}\nabla f_P(g^{-1}\cdot z),
\end{equation}
where $L_g:D \to D$ is the isometry given by left translation by $g \in G$.
Since $f_P$ is $\lsp{0}P$-invariant, 
\comment{eq:e2}
\begin{equation}\Label{eq:e2}
\nabla f_p(g^{-1} \cdot z)=(dL_{p^{-1}})_{pg^{-1} \cdot z}\nabla f_P(pg^{-1} \cdot z)\quad \text{for any $p\in \lsp{0}P$.}
\end{equation}
There exists a $p \in \lsp{0}P$ such that $p g^{-1}z=t \circ z$, where $t=\psi_P\left(\frac{f_P(g^{-1}z)}{f_P(z)}\right)$.  Note that 
\begin{equation}\Label{eq:taction}
f_P(t\circ z)=\left(\frac{f_P(g^{-1}z)}{f_P(z)}\right)f_P(z).
\end{equation}
Then by Lemma~\ref{lem:chainrule} and \eqref{eq:taction}
\begin{equation}\Label{eq:e3}
\nabla f_P(p g^{-1}\cdot z)=\left(\frac{f_P(g^{-1}z)}{f_P(z)}\right)d(t \circ \cdot)_z\nabla f_P(z).
\end{equation}
Combining \eqref{eq:e1}, \eqref{eq:e2}, and \eqref{eq:e3}, 
\begin{equation}
\nabla f_Q(z)=\frac{f_Q(z)}{f_P(z)}(dL_{g p^{-1}})_{p g^{-1} \cdot z}d(t \circ \cdot)_z\nabla f_P(z).
\end{equation}
Thus by Lemma~\ref{lem:geodesicisometry}, $f_P(z)\|\nabla f_Q(z)\|=f_Q(z)\|\nabla f_P(z)\|$.
\end{proof}

\begin{prop}\Label{prop:neggrad}
Let $P$ and $Q$ be distinct parabolic $\Q$-subgroups of $G$ and $z\in E(\{P,Q\})$.  Then the following are equivalent:
\begin{enumerate}
\item $z$ is a first contact point for $\{P,Q\}$.\Label{item:1}
\item $z$ is a critical point for $f_P|_{E(\{P,Q\})}$.\Label{item:2}
\item $\nabla  f_P(z)=-\nabla  f_Q(z)$.\Label{item:3}
\end{enumerate}
\end{prop}

\begin{proof}
Let $E=E(\{P,Q\})$.  By Proposition~\ref{prop:norminvariant}, $\|\nabla  f_P(z)\|=\|\nabla  f_Q(z)\|$ for $z\in E$.  It follows that 
\begin{equation*}
\nabla  f_P|_E(z)=\pr_{T_zE}(\nabla  f_P(z))=\frac{1}{2}(\nabla  f_P(z)+\nabla  f_Q(z)),
\end{equation*}
and hence \eqref{item:2} is equivalent to \eqref{item:3}.

It is clear that if $z$ a first contact for $\{P,Q\}$, $\nabla f_P(z)$ and $\nabla f_Q(z)$ point in opposite directions.  Thus Proposition~\ref{prop:norminvariant} implies that $\nabla f_P(z)=-\nabla f_Q(z)$ and hence \eqref{item:1} implies \eqref{item:2}.

Now suppose $\nabla f_P(z)=-\nabla f_Q(z)$ for some $z\in E(\{P,Q\})$.  By Proposition~\ref{prop:ckexplicit}, to show that $z$ is a first contact for $\{P,Q\}$, it suffices to show that $|c_k|=0$ for $k=0,\ldots N-1$.  Let $y= f_P(z)$ and fix $g\in G(\Q)$ such that $Q=\lsp{g}P$.  Then by Proposition~\ref{prop:fPexplicit}
\begin{equation*}
f_Q(z)=\frac{f_Q(g \cdot x_P)y}{\left(\sum_{k=0}^N |c_k(z)|^2 f_P(x^*)^{-2k} y^{2k}\right)^{1/N}} \quad \text{for all $z \in D$.}
\end{equation*}
Let $\frac{\partial}{\partial y}$ be the vector field on $D$ which is tangent to the flow of the geodesic action of $A_P$ such that $\frac{\partial}{\partial y}(f_P)=1$.  For $z\in E(\{P,Q\})$,
\begin{equation}\Label{eq:denom1}
\sum_{k=0}^N |c_k(z)|^2 f_P({x^*})^{-2k} y^{2k}=f_Q(g \cdot x_P)^N,
\end{equation}
and hence,
\begin{equation}\Label{eq:den}
 |c_N(z)|^2 f_P({x^*})^{-2N} y^{2N}=f_Q(g \cdot x_P)^N-\sum_{k=0}^{N-1} |c_k(z)|^2 f_P({x^*})^{-2k} y^{2k}.
\end{equation}
If $\nabla f_P(z)=-\nabla f_Q(z)$, then in particular, $\frac{\partial f_P}{\partial y}(z)=-\frac{\partial f_Q}{\partial y}(z)$. 
This implies that 
\begin{equation}
\begin{split}
1&=-\left(1-\frac{y}{Nf_Q(g \cdot x_P)^{N}}\sum_{k=0}^{N} 2k |c_k(z)|^2 f_P({x^*})^{-2k} y^{2k-1}\right)\\
&=-1+\sum_{k=0}^{N}\frac{2k|c_k(z)|^2}{Nf_Q(g \cdot x_P)^{N}f_P({x^*})^{2k}}y^{2k}.
\end{split}\Label{eq:partialy}
\end{equation}
Combining \eqref{eq:den} and \eqref{eq:partialy} gives
\begin{align}
1&=-1+2\left(1-\sum_{k=0}^{N-1}\left(\frac{N-k}{N}\right)\frac{|c_k(z)|^2}{f_Q(g \cdot x_P)^{N}f_P({x^*})^{2k}}y^{2k}\right).
\end{align}
Since $y>0$, $|c_k(z)|=0$ for $k=0, \ldots, N-1$.
\end{proof}

Propositions~\ref{prop:norminvariant} and \ref{prop:neggrad} immediately yield the following.
\begin{prop}\Label{prop:X00>X01}
Let $P$ and $Q$ be two distinct parabolics in $\Par$.  Then 
\[\langle \nabla  f_P,\nabla  f_P \rangle \geq |\langle \nabla  f_P,\nabla  f_Q\rangle|\]
on $E(\{P,Q\})$ with equality on the set of first contacts for $\{P,Q\}$.
\end{prop}
\end{section}
\begin{section}{The Spine}\Label{sec:spine}
For $z \in D$, denote by $M(z)$ the strongly admissible set $\I$ for which $z\in D'(\I)$. 

\begin{prop}\Label{prop:intersectingfiber}
Let $z$ be a point in $D$.  Then the set $M(z) \subset \Par$ is $\Stab{z}$-stable.   
\end{prop}
\begin{proof}
Let $P \in M(z)$ and $\gamma^{-1} \in \Stab{z}$.  Then $ f_P(z)= f_P(\gamma^{-1} \cdot z)=  f_{\lsp{\gamma} P}(z)$.  It follows that $\lsp{\gamma} P \in M(z)$. 
\end{proof}
\begin{thm}\Label{thm:retraction}
$D_0$ is a $\Gamma$-invariant deformation retract of $D$.
\end{thm}
\begin{proof}
Since $\Gamma$ acts on $\Par$ by conjugation, it also acts on subsets $\I$ of $\Par$ in the obvious way.  Suppose $z \in D_0$.  Then $z \in D'(\I)$ for some $\I$.  This means that 
\begin{enumerate}
\item $ f_{P}(z)= f_{Q}(z)\quad \text{for every pair $P, \ Q \in \I$ and}$
\item $ f_P(z)> f_R(z)\quad \text{for every $P \in \I$ and $R \in \Par \setminus \I$.}$
\end{enumerate}
The $\Gamma$-invariance of $\{f_P\}$ implies that 
\begin{enumerate}
\item $ f_{\lsp{\gamma}P}(\gamma \cdot z)= f_{\lsp{\gamma}Q}(\gamma \cdot z)\quad \text{for every pair $P, \ Q \in \I$ and}$
 \item $ f_{\lsp{\gamma}P}(\gamma \cdot z)> f_{\lsp{\gamma}R}(\gamma \cdot z)\quad \text{for every $P \in \I$ and $R \in \Par \setminus \I$.}$
\end{enumerate}
Thus $\gamma \cdot z \in D'({}^{\gamma}\I)$ for every $\gamma \in \Gamma$, and so $D_0$ is $\Gamma$-invariant.

Proposition~\ref{prop:tiling} and Definition~\ref{defn:spine} imply a decomposition of $D$ as a disjoint union
\begin{equation}
D=D_0\sqcup \coprod_{P \in \Par} D'(P)
\end{equation}

From Proposition~\ref{prop:fproperty},
\begin{equation*}
f_P(a \circ z)=a^\alpha f_P(z) \quad \text{for all $z\in D$ and $a\in A_P$.}
\end{equation*} 
Thus for $z\in D'(P)$, Proposition~\ref{prop:minf} implies that there exists an $a(z) \in A_P$ such that $a(z) \circ z \in D_0$.  Thus we can define a function $\mu:D \to \R_{>0}$ by $\mu(z)= a(z)^\alpha$, where it is understood that if $z\in D_0$ then $\mu(z)=1$.  

Recall that there is a natural isomorphism $\psi_P:\R_{>0} \to A_P$ defined in \eqref{eq:Raction}.  Define a family of maps $\tilde r_t:D \to D$ by
\begin{equation}
\tilde r_t(z)=\psi_P((1-t)+t \mu(z)) \circ z \quad \text{for $z \in D(P)$.}
\end{equation}
Note that this is well-defined because $\mu \equiv 1$ on $D_0$.  It is clear that $\tilde r_0=\id_D$ and $\tilde r_1(D) \subset D_0$.  The $\Gamma$-invariance of the exhaustion functions implies that
\begin{equation*}
\gamma \cdot \tilde r_t(z)=\tilde r_t(\gamma \cdot z) \quad \text{for every $\gamma \in \Gamma$ and $z\in D$.}
\end{equation*}

It is clear that $\tilde r$ is a continuous function of  $t$.  To see that it is a continuous function of $z$, it suffices to note that $\mu$ is a continuous function.  Thus $\tilde r$ gives a $\Gamma$-equivariant deformation retraction of $D$ onto $D_0$.  
\end{proof}

\begin{prop}\Label{prop:stabfiber}
For every $z\in D$, 
\begin{align*}
\Stab{z}&= \Stab{\tilde{r}_t(z)}\quad \text{for $t<1$ and}\\
\Stab{z}&\subseteq \Stab{\tilde{r}_1(z)}.  
\end{align*}
\end{prop}

\begin{proof}
Let $\gamma$ be an element of $\Stab{z}$.  Then $\gamma \cdot \tilde{r}_t(z) =\tilde r_t(\gamma \cdot z)=\tilde r_t(z)$.  Hence $\gamma\in \Stab{\tilde{r}_t(z)}$.  Notice that for each $z\notin D_0$, $c(t)=\tilde r_t(z), 0 \leq t \leq 1$ is a reparameterization of a geodesic, and $\Gamma$ acts by isometries.  Thus every $\gamma \in \Stab{z}$ fixes the geodesic through $\tilde r_1(z)$ and $z$.  In particular, $\Stab{z'}=\Stab{z}$ whenever $z=\tilde r_t(z')$ for some $t<1$.
\end{proof}

\begin{thm}\Label{thm:codim1and2}
If $\I \subset \Par$ is a subset of order two then $E(\I)$ is a contractible, smooth codimension $1$ submanifold of $D$.  If $\I$ has order three, then $E(\I)$ is a smooth codimension $2$ submanifold of $D$.
\end{thm}
\begin{proof}
Let $\I=\{P,Q\}$ and $E=E(\I)$.  Then $E$ is the zero set of $ f_P- f_Q$.  Thus to show smoothness, it suffices to show that $\nabla  f_P(z)-\nabla  f_Q(z) \neq 0$ for every $z \in E$. Suppose otherwise.  Consider the geodesic $c(t)$ with $c(0)=z$ and $c'(0)= \nabla  f_P(z)$.  For $t$ sufficiently large, $c(t)$ is contained in the degenerate tile $D'(P)$.  Similarly, since $c'(0)=\nabla  f_Q(z)$, for $t$ sufficiently large, $c(t)$ is in $D'(Q)$.  This contradicts the fact that $D'(P) \cap D'(Q)=\emptyset$.

To show contractibility, we will use the gradient flow of $ f_P|_E$ to contract E to the set of first contacts for $ f_P|_E$ and show that this set is contractible.  Proposition~\ref{prop:ckexplicit} implies that on this set, $ f_P(z)$ and $u(z)$ are constant, and $m(z)$ is free to range over $M_P(x)$.  In particular, the set of critical points is diffeomorphic to the symmetric space $D_{M_P}$ for $M_P$.  Since $D_{M_P}$ is contractible, this proves the result.

Now suppose $\I=\{P,Q,R\}$.  The set $E(\I)$ is defined by $h_1=0, \ h_2=0$, and $h_3=0$, where $h_1= f_P- f_Q, \ h_2= f_P- f_R$, and $h_3= f_Q- f_R$.  One must show that $\{\nabla h_1, \nabla h_2, \nabla h_3\}$ has constant rank on $E(\I)$.  Fix a point $z \in E(\I)$.  Note that $\nabla h_3=\nabla h_2 - \nabla h_1$, so the rank is less than 3.  All of $\nabla h_1,\ \nabla h_2,$ and $\nabla h_3$ are nonzero, so that the rank is either 1 or 2.  We prove the rank is not 1 by contradiction.  Suppose the rank is 1.  Then for each $z\in E(\I)$ there exist nonzero constants $d_1$ and $d_2$ such that 
\begin{align}\Label{eq:f1=cf3}
\nabla h_1(z)&=d_1\ \nabla h_3(z) \quad \text{and}\\
\Label{eq:f2=cf3}
\nabla h_2(z)&=d_2\ \nabla h_3(z).
\end{align}     

Using Proposition~\ref{prop:X00>X01} and taking the inner product of \eqref{eq:f1=cf3} with $\nabla  f_Q(z)$, it follows that $d_1<0$.  Similarly, by taking the inner product of \eqref{eq:f2=cf3} with $\nabla  f_R(z)$, it follows that $d_2>0$.  Taking the inner product of \eqref{eq:f1=cf3} with $\nabla  f_P(z)$ yields
\begin{equation*}
\ip{\nabla f_3(z)}{\nabla  f_P(z)} < 0.
\end{equation*}
Similarly, taking the inner product of \eqref{eq:f2=cf3} with $\nabla  f_P(z)$ yields
\begin{equation*}
\ip{\nabla f_3(z)}{\nabla  f_P(z)} > 0,
\end{equation*}
which gives the desired contradiction.  
\end{proof}

\begin{thm}\Label{thm:spine}
$D_0$ is a spine for $G$.  In particular, 
\begin{enumerate}
\item $D_0 \subset D$ is a $\Gamma$-equivariant deformation retract.
\item $\dim_\R(\Gamma \backslash D_0)=\vcd(\Gamma)$.
\item $\Gamma \backslash D_0$ is compact.
\end{enumerate}
\end{thm}

\begin{proof}
The first statement is Theorem~\ref{thm:retraction}.  By \cite{BS}, the virtual cohomological dimension of $G$ is given by 
\begin{align*}
\vcd(G)&=\dim_\R(D)-\rank_\Q(G) \\
&=\dim_\R(D)-1.
\end{align*}
Theorem~\ref{thm:codim1and2} implies that $\dim_\R(\Gamma \backslash D_0)=\vcd(G)$.  To see that $\Gamma \backslash D_0$ is compact, notice that $D_0 \subset D \subset \overline{D}$ is closed and $\Gamma \backslash \overline{D}$ is compact by \cite{BS}.
\end{proof}
\end{section}
\begin{section}{Separated groups}\Label{sec:separated}
One would like to know when the decomposition of $D_0$ given in Definition~\ref{defn:spine} is nice in some sense.  To this end, Theorem~\ref{thm:codim1and2} shows that the strongly admissible sets of order two correspond to open subsets of contractible, smooth, codimension 1 submanifolds of $D$, and that the strongly admissible sets of order three correspond to open subsets of smooth, codimension 2 submanifolds of $D$.  

The example of $\SL_2(\Z[\sqrt{2}])$, described in \cite{Brown} and recalled in Section~\ref{sec:examples} shows that it is too much to hope for that the strongly admissible sets always correspond smooth contractible sets.  Thus we try to give a criteria to check that ensures the strongly admissible sets correspond to smooth, contractible sets.  First we prove some geometric lemmas, and then we apply these lemmas to the gradient vector fields for the exhaustion functions to define a \emph{separated} condition that ensures that the strongly admissible sets correspond to smooth contractible sets. This condition is sufficient but not necessary.  We hope to define a refinement of $D_0$ to give a more satisfactory solution in a future paper.
\begin{subsection}{Geometric lemmas}
For this section, let $V$ be a finite dimensional real vector space with inner product $\ip{\cdot}{\cdot}$.  Fix a basis $\Delta=\{\alpha_1,\ldots\alpha_n\}$ of $V$ such that $\ip{\alpha_i}{\alpha_j} \leq 0$ for $i \neq j$.  Let $\hat \Delta =\{\beta_1,\ldots \beta_n\} \subset V$ be defined by $\ip{\beta_i}{\alpha_j}=\delta_{ij}$.  

\begin{defn}
Call a vector $v \in V$ \emph{dominant \textup (with respect to $\Delta$\textup )} if  $\ip{v}{\alpha_i}\geq 0$ for all $i$.  If these inequalities are strict for all $i$, then $v$ is \emph{strictly dominant}.  Call a vector $v \in V$ \emph{codominant \textup (with respect to $\Delta$\textup )} if  $\ip{v}{\beta_i}\geq 0$ for all $i$.  If these inequalities are strict for all $i$, then $v$ is \emph{strictly codominant}.     
\end{defn}

The dominant vectors in $V$ form a convex cone generated by $\hat \Delta$, while the codominant vectors form the dual cone generated by $\Delta$. 

We recall without proof a geometric lemma due to Langlands \cite{BoWa}.  
\begin{lem}\Label{lem:Langlands}
The dominant cone is contained within the codominant cone.  Equivalently, $\beta_i=\sum e_{ji}\alpha_j$ with $e_{ji}\geq 0$ for all $1\leq i,j,\leq n$.  Consequently, $\ip{\beta_i}{\beta_j}\geq 0$ for $1\leq i,j \leq n$.
\end{lem}

\begin{lem}\Label{lem:independent}
Let $\{v_0,\ldots v_k\}\subset V$ be a set of vectors such that $\ip{v_i}{v_j} <0$ for all $0 \leq i<j \leq k$.  Then the span of $\{v_0,\ldots v_k\}\subset V$ is either $k$ or $k+1$ dimensional.  In particular, any subset of $\{v_0,\ldots v_k\}$ of order $k$ is linearly independent.
\end{lem}

\begin{proof}
Let $W$ be the span of $S=\{v_0,\ldots, v_k\}$.  Let $\Delta \subseteq \{v_0,\ldots, v_k\}$ be a linearly independent subset.  Let $C$ be the dominant cone in $W$ with respect to $\Delta$.  Since $\ip{v_i}{v_j} <0$ for all $0 \leq i\leq j \leq k$, every vector in $S \setminus \Delta$ must lie in the anti-dominant cone $-C$.  By Lemma~\ref{lem:Langlands}, the vectors in $C$ have $-C$ have pairwise nonnegative inner product.  Hence $S \setminus \Delta$ has at most one element. \end{proof}

\begin{lem}\Label{lem:normals}
Let $\{v_0,\ldots, v_k\}\subset V$ be a set of vectors such that $\ip{v_i}{v_j} <0$ for all $0 \leq i<j \leq k$.  Then the span of $\{v_i-v_j\}_{0 \leq i,j\leq k}$ is a subspace of dimension $k$.  
\end{lem}

\begin{proof}
It is clear that the span of $\{v_i-v_j\}$ is equal to the span of $\{\eta_1,\ldots,\eta_k\}$ where $\eta_i=v_0-v_i$.  Thus it suffices to show that the $\eta_i$ are linearly independent.  Assume the contrary.  Then there exist $a_i \in \R$, not all zero, such that 
\begin{equation*}\sum_{i=1}^k a_i \eta_i=\sum_{i=1}^k a_i (v_0-v_i)=0.\end{equation*}
Then $\left(\sum_{i=1}^k a_i\right)v_0=\sum_{i=1}^k a_i v_i$.  By Lemma~\ref{lem:independent}, $v_1,\ldots, v_k$ are linearly independent so $\sum a_i\neq 0$.  Thus \[v_0=\sum_{i=1}^k\left(\frac{a_i}{\sum a_j}\right)v_i.\]  By Lemma~\ref{lem:Langlands}, $\frac{a_i}{\sum a_j} \leq 0$ for $1 \leq i \leq k$, which contradicts $\sum_{i=1}^k \frac{a_i}{\sum a_j}=1$. \end{proof}

\begin{lem}\Label{lem:boundary}
Let $\{v_0,\ldots, v_k\}\subset V$ be a set of linearly independent vectors such that $\ip{v_i}{v_j} \leq 0$ for all $0 \leq i<j \leq k$.  Let $u\in V$ be a vector such that $\ip{v_i}{u} \leq 0$ for all $0 \leq i \leq k$.  Then the orthogonal projection $\pr(v_0)$ of $v_0$ to the orthogonal complement in $V$ of $\operatorname{span}\{v_i-v_j\}_{0 \leq i,j \leq k}$ satisfies
\[\ip{\pr(v_0)}{v_i-u}>0 \quad \text{for all $0 \leq i \leq k$.}\]
\end{lem}

\begin{proof}
Let $\Delta=\{v_0,\ldots, v_k\}$ and $W \subset V$ the subspace spanned by $\Delta$.    Since the $v_i$ are linearly independent, $\pr(v_0) \neq 0$.  Furthermore, since $v_0-v_i$ is perpendicular to $\pr(v_0)$,
\begin{align*}
\ip{\pr(v_0)}{v_i}&=\ip{\pr(v_0)}{v_i+(v_0-v_i)}\\
&=\ip{\pr(v_0)}{v_0}\\
&=\|\pr(v_0)\|^2>0 \quad \text{for all $i$.}
\end{align*}
It follows that $\pr(v_0)$ lies in the dominant cone $C$ in $W$ with respect to $\Delta$ and has nonnegative inner product with every other vector in $C$ by Lemma~\ref{lem:Langlands}.  Since $\ip{v_i}{u} \leq 0$ for all $0 \leq i \leq k$, one can write $u$ as $u=\epsilon -u_C$, where $\epsilon$ is perpendicular to $W$ and $u_C$ lies in $C$.  Therefore 
\begin{align*}
\ip{\pr(v_0)}{v_i-u}&=\ip{\pr(v_0)}{v_i-(\epsilon -u_C)}\\
&=\|\pr(v_0)\|^2+\ip{\pr(v_0)}{u_C}\\
&>0 \quad \text{for all $i$.}\qedhere
\end{align*}
\end{proof}
\end{subsection}
\begin{subsection}{Consequences}
One can imagine the strongly admissible sets growing from their first contact as the parameter gets pushed lower.  Thus a retraction to the set of first contacts should be possible by ``going backward''.  More precisely, one would like to show that the negative gradient flow for $f_P$ (for any $P \in \I$), restricted to $D'(\I)$ defines a retraction onto the set of first contact points for $\I$.  Thus we need the restriction of the gradient flow to push each boundary piece into $D'(\I)$.    
\begin{lem}\Label{lem:generalsmooth}
Let $(M,g)$ be a Riemannian manifold and $f_0,\ldots,f_k$ be smooth functions defined on $M$.  Then the set 
\[\{p\in M \;| \; f_0(p)=\ldots =f_k(p) \} \cap \{p \in M \; |\; \ip{\nabla f_i}{\nabla f_j}_p <0 \text{ for $0\leq i<j \leq k$}\}\] is either empty or a smooth codimension $k$ submanifold of $M$.
\end{lem}

\begin{proof}
This is immediate from the Implicit Function Theorem and Lemma~\ref{lem:normals}.
\end{proof}

\begin{defn}\Label{defn:separated}
Given a strongly admissible set $\I$, let $\hat \I$ denote a minimal set of parabolics such that $D(\hat \I)=D(\I)$.  A strongly admissible set $\I$ is said to be \emph{separated} if $\ip{\nabla  f_P}{\nabla  f_Q} \leq 0$ on $D(\I)$ for every $P,~Q \in \hat \I$ and \emph{strictly separated} if the inequality is strict on $D(\I)$.  An algebraic group $\G$ is \emph{$($strictly$)$ separated} if its strongly admissible sets are (strictly) separated.
\end{defn}                                                                                                 
\begin{thm}\Label{thm:strictlyseparated}
Let $\G$ be a strictly separated group and $\I$ a strongly admissible set of order $k$.  Then $D'(\I)$ is a smooth codimension $k$ submanifold of $D$.  Furthermore, if $\overline{D'(\I)}$ is compact and there is a unique critical point for $f_\I$, then $D'(\I)$ is contractible.
\end{thm}

\begin{proof}
  The first statement follows from the definition of a strictly separated group and Lemma~\ref{lem:generalsmooth}.  

We claim that the gradient flow for $f_\I$ pushes the boundary of $D(\I)$ into the interior.  The boundary is a union of pieces of the form $D(\I \cup \{P\})$ for some rational parabolic subgroup $P$.  This piece flows into the interior of $D(\I)$ if 
\begin{equation}\Label{eq:ip}
\ip{\nabla f_\I}{\nabla f_\I-\nabla f_P}>0.
\end{equation}
Write $\I=\{Q_0, \ldots , Q_k\}$ and apply Lemma~\ref{lem:boundary} with $v_i=\nabla f_{Q_i}$ for $i=0,\ldots , k$ and $u=\nabla f_P$ to show that \eqref{eq:ip} is satisfied.  Thus the gradient flow of $f_\I$ will define a deformation retraction of $D(\I)$ onto the critical point.
\end{proof}
\end{subsection}
\end{section}

\begin{section}{Examples}\Label{sec:examples}
We apply the methods of Sections~\ref{sec:definitions}-\ref{sec:properties} to briefly describe the spines in some examples and compare our method with other known methods. 
\begin{subsection}{Bianchi groups}
The arithmetic group $\SL_2(\OOO_k)$ for $k$ a imaginary quadratic field viewed as a subgroup of $\SL_2(\C)$ acts on hyperbolic 3-space $\mathbb{H}^3$.  $\mathbb{H}^3$ is a linear symmetric space, and hence the existence of a spine (the well-rounded retract) is guaranteed by \cite{A2}.  Mendoza \cite{Men} and Vogtmann \cite{Vog} use a notion of distance to a cusp due to Siegel \cite{Sie} to explicitly compute a spine.  Ash \cite{A} shows that the spines they compute are the same as the well-rounded retract.  The distance functions are a power of the exhaustion functions that we define in \ref{defn:exhaustion}, and thus the spines also coincide with our spines.   
\end{subsection}
\begin{subsection}{Hilbert modular group}\label{subsec:hilbert}
The arithmetic group $\SL_2(\OOO_k)$, for $k$ a real quadratic field, acts on the product of 2 upper half-planes $\HH\times\HH$.  $\HH\times\HH$ is a linear symmetric space, and hence the existence of a spine (the well-rounded retract) is guaranteed by \cite{A2}.  Brownstein again uses Siegel's notion of distance to a cusp to explicitly compute the spine in \cite{Brown} for $k=\Q(\sqrt{5})$.  He also conjectures a spine for $k=\Q(\sqrt{2})$.  One can easily show that these spines are different from those obtained by Ash.  On the other hand, the distance functions considered by Brownstein are a power of the exhaustion functions that we define in \ref{defn:exhaustion}, and thus Brownstein's spines also coincide with ours.
\begin{figure}
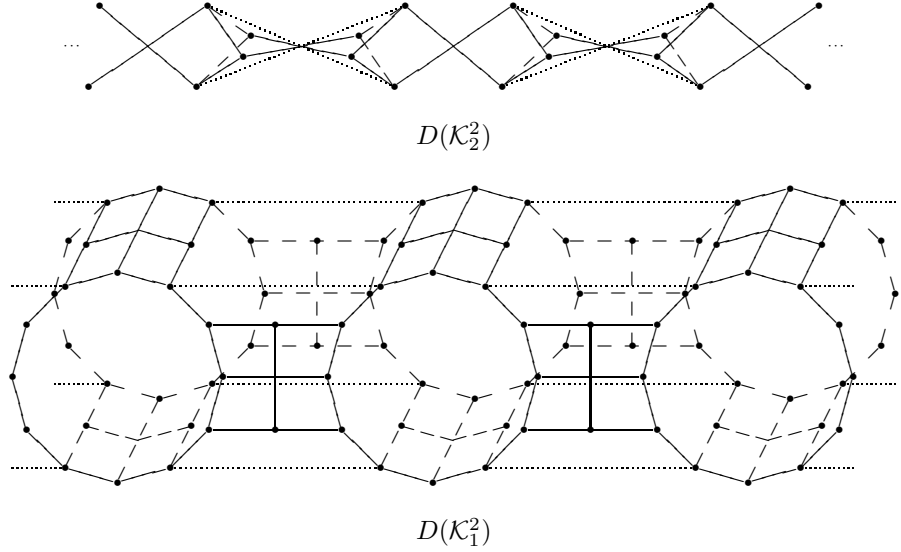

\[\begin{array}{c}
{\def\objectstyle{\scriptscriptstyle}
\xy /r0.4in/:="b", +(2,0)="c",+(-4,0)="a",+(-2,0)="z",+(8,0)="y",+(1,0)="u",+(-10,0)="v",
"a",{\xypolygon4"A"{~:{(1,0):(.07,.2)::}~>{}{\bullet}}},
"b",{\xypolygon4"B"{~:{(1,0):(.1,.75)::}~>{}{\bullet}}},
"c",{\xypolygon4"C"{~:{(1,0):(.07,.2)::}~>{}{\bullet}}},
"z",{\xypolygon4"Z"{~:{(1,0):(.1,.75)::}~>{}{\bullet}}},
"y",{\xypolygon4"Y"{~:{(1,0):(.1,.75)::}~>{}{\bullet}}},
"u"*{\cdots},"v"*{\cdots},
"B1";"B3"**@{-},"B2";"B4"**@{-},
"A1";"A3"**@{-},"A2";"A4"**@{-},
"C1";"C3"**@{-},"C2";"C4"**@{-},
"Z1";"Z3"**@{-},"Z2";"Z4"**@{-},
"Y1";"Y3"**@{-},"Y2";"Y4"**@{-},
"B1";"C2"**@{--},"B4";"C2"**@{--},"B1";"C3"**@{-},"B4";"C3"**@{-},
"B2";"A4"**@{-},"B3";"A4"**@{-},"B2";"A1"**@{--},"B3";"A1"**@{--},
"B2";"Z4"**@{.},"B3";"Z1"**@{.},"B1";"Y3"**@{.},"B4";"Y2"**@{.},
"Y2";"C4"**@{-},"Y3";"C4"**@{-},"Y2";"C1"**@{--},"Y3";"C1"**@{--},
"Z1";"A2"**@{--},"Z4";"A2"**@{--},"Z1";"A3"**@{-},"Z4";"A3"**@{-},
\endxy
}\\
\\
D(\K^2_2)\\
\\
{\def\objectstyle{\scriptscriptstyle}
\xy /r.55in/:="bd", +(3,0)="cd",+(-6,0)="ad",
"ad",{\xypolygon12"AD"{~={0}{\bullet}}},
"bd",{\xypolygon12"BD"{~={0}{\bullet}}},
"cd",{\xypolygon12"CD"{~={0}{\bullet}}},
"bd",+(0.4,.8)="bu",{\xypolygon12"BU"{~={0}~>{}{\bullet}}},
"ad",+(0.4,.8)="au",{\xypolygon12"AU"{~={0}~>{}{\bullet}}},
"cd",+(0.4,.8)="cu",{\xypolygon12"CU"{~={0}~>{}{\bullet}}},
"bd",+(0.2,.4)="bm",{\xypolygon12"BM"{~={0}~>{}}},
"ad",+(0.2,.4)="am",{\xypolygon12"AM"{~={0}~>{}}},
"cd",+(0.2,.4)="cm",{\xypolygon12"CM"{~={0}~>{}}},
"AD3";"AU3"**@{-},"AD4";"AU4"**@{-},"AD5";"AU5"**@{-},
"AM3",*{\bullet},"AM5",*{\bullet},
"AM3";"AM4"**@{-},"AM5";"AM4"**@{-},
"AU3";"AU4"**@{-},"AU5";"AU4"**@{-},
"BD3";"BU3"**@{-},"BD4";"BU4"**@{-},"BD5";"BU5"**@{-},
"BM3",*{\bullet},"BM5",*{\bullet},
"BM3";"BM4"**@{-},"BM5";"BM4"**@{-},
"BU3";"BU4"**@{-},"BU5";"BU4"**@{-},
"CD3";"CU3"**@{-},"CD4";"CU4"**@{-},"CD5";"CU5"**@{-},
"CM3",*{\bullet},"CM5",*{\bullet},
"CM3";"CM4"**@{-},"CM5";"CM4"**@{-},
"CU3";"CU4"**@{-},"CU5";"CU4"**@{-},
"AD9";"AU9"**@{--},"AD10";"AU10"**@{--},"AD11";"AU11"**@{--},
"AM9",*{\bullet},"AM11",*{\bullet},
"AM9";"AM10"**@{--},"AM11";"AM10"**@{--},
"AU9";"AU10"**@{--},"AU11";"AU10"**@{--},
"BD9";"BU9"**@{--},"BD10";"BU10"**@{--},"BD11";"BU11"**@{--},
"BM9",*{\bullet},"BM11",*{\bullet},
"BM9";"BM10"**@{--},"BM11";"BM10"**@{--},
"BU9";"BU10"**@{--},"BU11";"BU10"**@{--},
"CD9";"CU9"**@{--},"CD10";"CU10"**@{--},"CD11";"CU11"**@{--},
"CM9",*{\bullet},"CM11",*{\bullet},
"CM9";"CM10"**@{--},"CM11";"CM10"**@{--},
"CU9";"CU10"**@{--},"CU11";"CU10"**@{--},
"AD12";"BD8"**@{-}?="abdb"*{\bullet},"AD1";"BD7"**@{-},"AD2";"BD6"**@{-}?="abdt"*{\bullet},
"AU12";"BU8"**@{--},"AU1";"BU7"**@{--},"AU2";"BU6"**@{--},
"AU12";"AU1"**@{--},"AU2";"AU1"**@{--},
"BU8";"BU7"**@{--},"BU6";"BU7"**@{--},
"abdb"+(.4,.8)="abub"*{\bullet},"abdt"+(.4,.8)="abut"*{\bullet},
"abdb";"abdt"**@{-},"abub";"abut"**@{--},
"BD12";"CD8"**@{-}?="bcdb"*{\bullet},"BD1";"CD7"**@{-},"BD2";"CD6"**@{-}?="bcdt"*{\bullet},
"BU12";"CU8"**@{--},"BU1";"CU7"**@{--},"BU2";"CU6"**@{--},
"BU12";"BU1"**@{--},"BU2";"BU1"**@{--},
"CU8";"CU7"**@{--},"CU6";"CU7"**@{--},
"bcdb"+(.4,.8)="bcub"*{\bullet},"bcdt"+(.4,.8)="bcut"*{\bullet},
"bcdb";"bcdt"**@{-},"bcub";"bcut"**@{--},
"AU5";"AU6"**@{--},"AU6";"AU7"**@{--},"AU7";"AU8"**@{--},"AU9";"AU8"**@{--},
"AU3";"AU2"**@{--},"AU12";"AU11"**@{--},
"BU5";"BU6"**@{--},"BU9";"BU8"**@{--},
"BU3";"BU2"**@{--},"BU11";"BU12"**@{--},
"CU5";"CU6"**@{--},"CU9";"CU8"**@{--},
"CU3";"CU2"**@{--},"CU1";"CU2"**@{--},"CU1";"CU12"**@{--},"CU11";"CU12"**@{--},
"AU3";"BU5"**@{.},"BU3";"CU5"**@{.},
"AD3";"BD5"**@{.},"BD3";"CD5"**@{.},
"AD11";"BD9"**@{.},"BD11";"CD9"**@{.},
"AU11";"BU9"**@{.},"BU11";"CU9"**@{.},
"AD9"+(-0.5,0);"AD9"**@{.},"AD5"+(-0.5,0);"AD5"**@{.},
"AU9"+(-0.5,0);"AU9"**@{.},"AU5"+(-0.5,0);"AU5"**@{.},
"CD11"+(0.5,0);"CD11"**@{.},"CD3"+(0.5,0);"CD3"**@{.},
"CU11"+(0.5,0);"CU11"**@{.},"CU3"+(0.5,0);"CU3"**@{.},
\endxy}\\
\\
D(\K^2_1)
\end{array}\]
\caption{The two $\SL_2(\Z[\sqrt{2}])$-conjugacy classes of strongly admissible sets of order two are both infinite in length, and the boundary pieces can be identified with certain configurations of vectors in $\Z[\sqrt{2}]^2$.  The boundary of $D(\K^2_2)$ consists of pieces which are the union of two quadrilaterals joined at a vertex.  The boundary of $D(\K^2_1)$ consists of dodecagons and the union of quadrilaterals described above.}
\label{fig:hilbertround}
\end{figure}

We consider the example of $\Q(\sqrt{2})$ in some detail.  Fix the following sets $\LL^n_m$ of cusps in $k \cup \{\infty\}\simeq \OOO^*\backslash \OOO^2$.
\begin{gather*}
\LL^2_1=\{\infty,0\} \qquad\LL^2_2=\{\infty,\frac{1}{\sqrt{2}}\}\\
\LL^3_1=\{\infty,0,1\}\qquad\LL^3_2=\{\infty,0,\frac{1}{\sqrt{2}}\}\\
\LL^4_1=\{\infty,0,1,\frac{1}{\sqrt{2}}\}\qquad\LL^4_2=\{\infty,0,\sqrt{2},\frac{1}{\sqrt{2}}\}\\
\LL^5=\{\infty,0,1,\sqrt{2},\frac{1}{\sqrt{2}}\}
\end{gather*}
Let $\K^n_m$ denote the set of rational parabolics associated to $\LL^n_m$.  Then $D(\K^n_m) \subset \HH\times \HH$ is defined by the equality of $n$ exhaustion functions and has codimension equal to $n-1$.

These subsets of cusps correspond to geometrically different pieces of the spine.  However, a given type may contain more than one $\Gamma$-conjugacy class of strongly admissible set.  $D(\K^5)$ is a vertex.  $D(\K^4_1)$ is a line segment, and $D(\K^4_2)$ is the union of 2 transverse line segments.  $D(\K^3_1)$ is a dodecagon, and $D(\K^3_2)$ is the union of two quadrilaterals joined at a vertex.  Figure~\ref{fig:hilbertround} shows that $D(\K^2_1)$ is an infinite tube with boundary faces of types $D(\K^3_1)$ and $D(\K^3_2)$, and $D(\K^2_2)$ is homeomorphic to an infinite string of 3-cells with boundary faces of type $D(\K^3_2)$. The strongly admissible sets of order greater than two have finite stabilizers.  On the other hand, stabilizers of the strongly admissible sets of order two contain a subgroup of finite index that is isomorphic to the group of units $\OOO^*_k$, acting by translation along the infinite direction.  However, it is clear that the sets may be refined so that the spine has the structure of a regular cell complex, on which $\Gamma$ acts with finite stabilizers.  

The incidence table is given in Table~\ref{tab:hilbertincidence}, where the entry below the diagonal means that each column cell has that many row cells in its boundary, and the entry above the diagonal means the column cell appears in the boundary of this many row cells.
\begin{table}
\caption{Incidence types}
\label{tab:hilbertincidence}
\begin{tabular}{|c|cc|cc|cc|c|}
\hline
 & $\K^2_1$&$\K^2_2$&$\K^3_1$&$\K^3_2$&$\K^4_1$&$\K^4_2$&$\K^5$\\\hline
$\K^2_1$&$\ast$&$\ast$&3 &2 &5 & 4  & 8\\
$\K^2_2$&$\ast$&$\ast$&0 &1 &1 & 2  & 2\\\hline
$\K^3_1$&$\infty$ & 0&$\ast$&$\ast$& 2& 0  & 4\\
$\K^3_2$& $\infty$& $\infty$&$\ast$&$\ast$& 2&4   &6 \\\hline
$\K^4_1$&$\infty$ & $\infty$&12 &4 &$\ast$&$\ast$&4  \\
$\K^4_2$&$\infty$ & $\infty$& 0& 1&$\ast$&$\ast$&1  \\\hline
$\K^5$&$\infty$ & $\infty$& 12& 6& 2& 4&$\ast$\\
\hline
\end{tabular}
\end{table}
\end{subsection}
\begin{subsection}{Picard modular group}
The symmetric space of $\SU(2,1)$ is a complex 2-ball.  It is not a linear symmetric space, and hence the well-rounded retract does not apply for this group.  Our method of exhaustion functions applied to this group exhibits the second known example of a spine in the non-linear case (after $\Sp_4(\Z)$ by MacPherson and McConnell in \cite{MM}).  Details about this spine, including stabilizer information and cohomology computations can be found in \cite{Yaspicard}.  

\begin{figure}
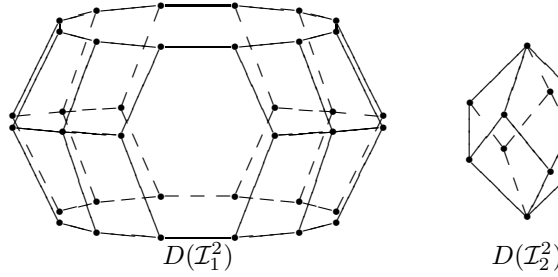

\[\begin{array}{c@{\hspace{0.4in}}c}
{\def\objectstyle{\scriptscriptstyle}
\xy /r0.5in/:="a", +(0,1)="b",+(0,1)="c","a",
{\xypolygon12"A"{~:{(1.5,0):(0,.15)::}~>{}{\bullet}}},
"b",
{\xypolygon12"B"{~:{(2,0):(0,.075)::}~>{}}},
"b",
{\xypolygon8"D"{~:{(2.1,0):(0,.075)::}~>{}{\bullet}}},
"c",
{\xypolygon12"C"{~:{(1.5,0):(0,.15)::}{\bullet}}},
"B2"*{\bullet},"B5"*{\bullet},"B8"*{\bullet},"B11"*{\bullet},
"A1";"D1"**@{--},"A2";"B2"**@{--},"A3";"D2"**@{--},"A4";"D3"**@{--},
"A5";"B5"**@{--},"A6";"D4"**@{--},"A7";"D5"**@{-},"A8";"B8"**@{-},
"A9";"D6"**@{-},"A10";"D7"**@{-},"A11";"B11"**@{-},"A12";"D8"**@{-},
"C1";"D1"**@{-},"C2";"B2"**@{--},"C3";"D2"**@{--},"C4";"D3"**@{--},
"C5";"B5"**@{--},"C6";"D4"**@{-},"C7";"D5"**@{-},"C8";"B8"**@{-},
"C9";"D6"**@{-},"C10";"D7"**@{-},"C11";"B11"**@{-},"C12";"D8"**@{-},
"D5";"B8"**@{-},"D6";"B8"**@{-},"D7";"B11"**@{-},"B11";"D8"**@{-},
"D5";"B8"**@{-},"D6";"B8"**@{-},"D7";"B11"**@{-},"B11";"D8"**@{-},
"D4";"B5"**@{--},"B5";"D3"**@{--},"D2";"B2"**@{--},"B2";"D1"**@{--},
"A7";"A6"**@{--},"A6";"A5"**@{--},"A5";"A4"**@{--},"A3";"A4"**@{--},"A3";"A2"**@{--},"A2";"A1"**@{--},"A1";"A12"**@{--},
"A7";"A8"**@{-},"A9";"A8"**@{-},"A9";"A10"**@{-},"A10";"A11"**@{-},"A11";"A12"**@{-}
\endxy}&
{
\def\objectstyle{\scriptscriptstyle}
\xy /r0.3in/:="a", +(0,1)="b",+(0,1)="c",+(0,1)="d",
"a"*{\bullet},"d"*{\bullet},
"b",{\xypolygon4"B"{~={45}~:{(1,-0.14):(0.4,.2)::}~>{}{\bullet}}},
"c",{\xypolygon4"C"{~={0}~:{(1,0):(.4,.2)::}~>{}{\bullet}}},
"C1";"B1"**@{-},"C1";"B4"**@{-},
"C2";"B2"**@{--},"C2";"B1"**@{--},
"C3";"B3"**@{-},"C3";"B2"**@{--},
"C4";"B4"**@{-},"C4";"B3"**@{-},
"B1";"a"**@{-},"B2";"a"**@{--},"B3";"a"**@{-},"B4";"a"**@{-},
"C1";"d"**@{-},"C2";"d"**@{--},"C3";"d"**@{-},"C4";"d"**@{-}
\endxy}\\
D(\I^2_1)&D(\I^2_2)
\end{array}\]
\caption{The two $\SU(2,1;\Z[i])$-conjugacy classes of strongly admissible sets of order two.  $D(\I^2_1)$ is homeomorphic to a polytope with dodecagon, hexagon and quadrilateral faces, while $D(\I^2_2)$ has only quadrilateral faces.  The strongly admissible sets are parameterized by certain configurations of vectors in $\Z[i]^3$.}
\label{fig:picardround}
\end{figure}

We can identify the cusps with isotropic vectors in $\Z[i]^3$ modulo scaling by $\Z[i]^*$.  Fix the following sets $\J^n_m$ of cusps.
\begin{gather*}
\J^2_1=\left\{\vect{1}{0}{0},\vect{0}{0}{1}\right\} \qquad\J^2_2=\left\{\vect{1}{0}{0},\vect{i}{1+i}{1+i}\right\}\\
\J^3_1=\J^2_1 \cup \left\{\vect{1}{0}{1}\right\}\qquad\J^3_2=\J^2_1 \cup\left\{\vect{i}{1+i}{1}\right\}\qquad \J^3_3=\J^2_1  \cup\left\{\vect{1+i}{1+i}{1}\right\}\\
\J^4_1=\J^3_1\cup\left\{\vsigma \right\} \qquad \J^4_2=\J^3_3\cup \left\{\vect{-1}{-1+i}{1+i}\right\}\\
\J^5=\J^4_1 \cup\left\{\vect{1+i}{1+i}{1}\right\}\\
\J^8=\left\{\ve,\vw,\begin{pmatrix} -1\\1+i\\1+i \end{pmatrix},\vect{-1+i}{1+i}{1},\vect{1+i}{1-i}{1},\begin{pmatrix} i\\1+i\\1+i \end{pmatrix},\vect{2i}{2}{1},\vect{i}{2}{2}\right\}.
\end{gather*}
Let $\I^n_m$ denote the set of rational parabolics associated to $\J^n_m$.  

These subsets of cusps correspond to geometrically different pieces of the spine.  Each class contains exactly one $\Gamma$-conjugacy class of strongly admissible set.  $D(\I^8)$ is a vertex, and vertices that occur where four quadrilaterals come together are all $\Gamma$-conjugate to $D(\I^8)$.  The other vertices are conjugate to $D(\I^5)$.  The $\Gamma$-translates of $D(\I^4_1)$ are the edges that do not terminate in a vertex of type $D(\I^8)$.  The other edges have type $D(\I^4_2)$.  $D(\I^3_1)$ is a dodecagon, and $D(\I^3_3)$ is a hexagon.  Figure~\ref{fig:picardround} shows $D(\I^2_1)$ and $D(\I^2_2)$. They are both homeomorphic to polytopes.

The incidence table is given in Table~\ref{tab:picardincidence}, where as in \ref{subsec:hilbert} the entry below the diagonal means that each column cell has that many row cells in its boundary, and the entry above the diagonal means the column cell appears in the boundary of this many row cells.
\begin{table}

\caption{Incidence Types}
\label{tab:picardincidence}
\begin{tabular}{|c|cc|ccc|cc|cc|}
\hline
 & $\I^2_1$&$\I^2_2$&$\I^3_1$&$\I^3_2$&$\I^3_3$&$\I^4_1$&$\I^4_2$&$\I^5$&$\I^8$\\\hline
$\I^2_1$&$\ast$&$\ast$&3&3&2&5&4&8&16\\
$\I^2_2$&$\ast$&$\ast$&0&0&1&1&2&2&8\\\hline
$\I^3_1$&2&0&$\ast$&$\ast$&$\ast$&1&0&2&0\\
$\I^3_2$&4&0&$\ast$&$\ast$&$\ast$&1&0&2&0\\
$\I^3_3$&12&8&$\ast$&$\ast$&$\ast$&2&4&6&32\\\hline
$\I^4_1$&40&8&12&6&2&$\ast$&$\ast$&4&0\\
$\I^4_2$&16&8&0&0&2&$\ast$&$\ast$&1&16\\\hline
$\I^5$&32&8&12&6&3&2&1&$\ast$&$\ast$\\
$\I^8$&4&2&0&0&1&0&1&$\ast$&$\ast$\\
\hline
\end{tabular}
\end{table}
\end{subsection}
\end{section}
\bibliography{../references}    
\end{document}